\documentclass[12pt]{article}

\usepackage{amsmath,amsfonts,amsthm,amssymb}
\usepackage{fullpage}

\theoremstyle{definition}
\newtheorem{definition}{Definition}
\theoremstyle{remark}

\newtheorem{example}[definition]{Example}

\newtheoremstyle{mytheorem}{0.5cm}{0.2cm}{\slshape}{ }{\bfseries}{.}{ }{}
\theoremstyle{mytheorem}
\newtheorem{theorem}[definition]{Theorem}
\newtheorem{prop}[definition]{Proposition}

\newtheorem{cor}[definition]{Corollary}

\newcommand{\stsets}[1]{\mathbb{#1}}
\newcommand{\R}{\stsets{R}}
\newcommand{\Rd}{\R^d}
\renewcommand{\SS}{{\stsets{S}_+}}
\newcommand{\Sphere}[1][d-1]{\stsets{S}^{#1}}
\newcommand{\EE}{\stsets{E}}
\newcommand{\LL}{\stsets{L}}
\newcommand{\II}{\stsets{I}}

\newcommand{\Prob}[1]{\mathbf{P}\{#1\}}
\DeclareMathOperator{\E}{{\bf E}}

\DeclareMathOperator{\conv}{conv}
\DeclareMathOperator{\one}{{ 1\hspace*{-0.55ex}I}}

\newcommand{\simplex}[1]{\varDelta_{#1}}
\newcommand{\mes}{V}

\newcommand{\comp}{\mathbf{c}}

\renewcommand{\phi}{\varphi}
\newcommand{\ti}{\to\infty}

\newcommand{\dsim}{\overset{\mathrm{d}}{\sim}}

 \newlength{\querylen}
 \usepackage{fancybox}

\numberwithin{equation}{section}
\numberwithin{definition}{section}

\begin{document}
\bibliographystyle{plain}

\title{Convex geometry of max-stable distributions}
\author{\textsc{Ilya Molchanov}\\
  \normalsize
  Department of Mathematical Statistics and Actuarial Science,\\
  \normalsize
  University of Bern, Sidlerstrasse 5, CH-3012 Bern, Switzerland\\
  \normalsize
  E-mail: ilya@stat.unibe.ch
}
\date{}
\maketitle

\begin{abstract}
  \noindent 
  It is shown that max-stable random vectors in $[0,\infty)^d$ with
  unit Fr\'echet marginals are in one to one correspondence with
  convex sets $K$ in $[0,\infty)^d$ called max-zonoids. The
  max-zonoids can be characterised as sets obtained as limits of
  Minkowski sums of cross-polytopes or, alternatively, as the
  selection expectation of a random cross-polytope whose distribution
  is controlled by the spectral measure of the max-stable random
  vector.  Furthermore, the cumulative distribution function
  $\Prob{\xi\leq x}$ of a max-stable random vector $\xi$ with unit
  Fr\'echet marginals is determined by the norm of the inverse to $x$,
  where all possible norms are given by the support functions of
  (normalised) max-zonoids. As an application, geometrical
  interpretations of a number of well-known concepts from the theory
  of multivariate extreme values and copulas are provided.

  \medskip

\end{abstract}

\newpage

\section{Introduction}
\label{sec:introduction}

A random vector $\xi$ in $\Rd$ is said to have a \emph{max-stable}
distribution if, for every $n\geq2$, the cooordinatewise maximum of
$n$ i.i.d. copies of $\xi$ coincides in distribution with an affine
transform of $\xi$, i.e.
\begin{equation}
  \label{eq:ms}
  \xi^{(1)}\vee\cdots\vee\xi^{(n)} \dsim a_n\xi+b_n
\end{equation}
for $a_n>0$ and $b_n\in\Rd$. If (\ref{eq:ms}) holds with $b_n=0$ for
all $n$, then $\xi$ is called \emph{strictly} max-stable, see, e.g.,
\cite{beir:goeg:seg:04,kot:nad00,res87}.

Since every max-stable random vector $\xi$ is infinitely divisible
with respect to coordinatewise maximum, its cumulative distribution
function satisfies
\begin{equation}
  \label{eq:cdf}
  F(x)=\Prob{\xi\leq x}=
  \begin{cases}
    \exp\{-\mu([-\infty,x]^\comp)\}\,, & x\geq a\,,\\
    0\,, & \text{otherwise}\,,
  \end{cases}
  \quad x\in\Rd\,,
\end{equation}
where $a\in[-\infty,\infty)^d$, the superscript $\comp$ denotes the
complement and $\mu$ is a measure on $[a,\infty]\setminus\{a\}$ called
the \emph{exponent} measure of $\xi$, see \cite[Prop.~5.8]{res87}.
Note that all inequalities and segments (intervals) for vectors are
understood coordinatewise.

Representation (\ref{eq:cdf}) shows that the cumulative distribution
function of $\xi$ can be represented as the exponential
$F(x)=e^{-\nu(x)}$ of another function $\nu$. If $\xi$ is strictly
max-stable and $a=0$, then $\nu$ is homogeneous, i.e.
$\nu(sx)=s^{-\alpha} \nu(x)$ for all $s>0$ and some $\alpha>0$.  This
fact can be also derived from general results concerning
semigroup-valued random elements \cite{dav:mol:zuy05}. If $\alpha=1$,
an example of such function $\nu(x)$ is provided by $\nu(x)=\|x^*\|$,
i.e. a norm of $x^*=(x_1^{-1},\dots,x_d^{-1})$ for
$x=(x_1,\dots,x_d)\in[0,\infty)^d$. One of the main aims of this paper
is to show that this is the \emph{only} possibility and to
characterise all norms that give rise to strictly max-stable
distributions with $\alpha=1$. 

Every norm is homogeneous and sublinear. It is known
\cite[Th.~1.7.1]{schn} that each bounded homogeneous and sublinear
function $g$ on $\Rd$ can be described as the \emph{support function}
of a certain convex compact set $K$, i.e.
\begin{displaymath}
  g(x)=h(K,x)=\sup\{\langle x,y\rangle:\; y\in K\}\,,
\end{displaymath}
where $\langle x,y\rangle$ is the scalar product of $x$ and $y$. In
Section~\ref{sec:convex-sets-related} we show that every standardised
strictly max-stable distribution with $\alpha=1$ is associated with
the unique compact convex set $K\subset[0,\infty)^d$ called the
dependency set. The dependency sets are suitably rescaled sets from
the family of sets called max-zonoids.  While classical zonoids appear
as limits for the sums of segments \cite[Sec.~3.5]{schn}, max-zonoids
are limits of the sums of cross-polytopes. The contributions of
particular cross-polytopes to this sum are controlled by the spectral
measure of the max-stable random vector.  It is shown that not every
convex compact set for $d\geq3$ corresponds to a strictly max-stable
distribution, while if $d=2$, then the family of dependency sets is
the family of all ``standardised'' convex sets, see also \cite{falk06}
for the treatment of the bivariate case. This, in particular, shows a
substantial difference between possible dependency structures for
bivariate extremes on one hand and multivariate extremes in dimensions
three and more on the other hand.

The geometrical interpretation of max-stable distributions opens a
possibility to use tools from convex geometry in the framework of the
theory of extreme values. For instance, the polar sets to the
dependency set $K$ appear as multivariate quantiles of the
corresponding max-stable random vector, i.e. the level sets of its
cumulative distribution function.  In the other direction, some useful
families of extreme values distributions may be used to construct new
norms in $\R^d$ which acquire an explicit probabilistic interpretation.
The norms corresponding to max-stable distributions are considered in
Section~\ref{sec:norm-associated-with}.

Section~\ref{sec:expon-meas-spectr} deals with relationships between
spectral measures of max-stable laws and geometric properties of the
corresponding dependency set. In Section~\ref{sec:copulas-association}
it is shown that a number of dependency concepts for max-stable random
vectors can be expressed using geometric functionals of the dependency
set and its polar. Here also relationships to copulas are considered.
It is shown that max-zonoids are only those convex sets whose support
functions generate multivariate extreme value copulas.

It is well known that $Z$ is (classical) zonoid if and only if
$e^{-h(Z,x)}$ is positive definite, see \cite[p.~194]{schn}. In
Section~\ref{sec:compl-altern-extr} we establish a similar result for
the positive definiteness of the exponential with respect to the
coordinatewise maximum operation in case $Z$ is a max-zonoid. 

Section~\ref{sec:oper-with-conv} describes some relationships between
operations with convex sets and operations with max-stable random
vectors. Finally, Section~\ref{sec:infin-dimens-case} briefly
mentions an infinite-dimensional extension for max-stable sample
continuous random functions.

\section{Dependency sets and max-zonoids}
\label{sec:convex-sets-related}

Let $\xi$ be a max-stable random vector with non-degenerate marginals.
By an affine transformation it is possible to standardise the
marginals of $\xi$, so that $\xi$ has $\Phi_\alpha$ (Fr\'echet
distributed) marginals, where
\begin{displaymath}
  \Phi_\alpha(x)=
  \begin{cases}
    0, & x<0\,,\\
    e^{-x^{-\alpha}}, & x\geq0\,,
  \end{cases}
  \quad \alpha>0\,,
\end{displaymath}
or $\Psi_\alpha$ (Weibull or negative exponential distributed)
marginals, i.e.
\begin{displaymath}
  \Psi_\alpha(x)=
  \begin{cases}
    e^{-(-x)^\alpha}, & x<0\,,\\
    1, & x\geq 0\,,
  \end{cases}
  \quad \alpha>0\,,
\end{displaymath}
or $\Lambda$ (Gumbel or double exponentially distributed) marginals, i.e.
\begin{displaymath}
  \Lambda(x)=\exp\{-e^{-x}\}\,, x\in\R\,.
\end{displaymath}
By using (possibly non-linear) monotonic transformations applied to
the individual coordinates it is possible to assume that all marginals
are $\Phi_1$, see \cite[Prop.~5.10]{res87} and
\cite[Sec.~8.2.2]{beir:goeg:seg:04}. In this case we say that $\xi$
has \emph{unit Fr\'echet marginals} or has a \emph{simple} max-stable
distribution, see also \cite{ein:haan:sin97}. Sometimes we say that
$\xi=(\xi_1,\dots,\xi_d)$ has a \emph{semi-simple} max-stable
distribution if its rescaled version $(c_1\xi_1,\dots,c_d\xi_d)$ has a
simple max-stable distribution for some $c_1,\dots,c_d>0$.

If $\xi$ has a simple max-stable distribution, then
\cite[Prop.~5.11]{res87} implies that the exponent in
(\ref{eq:cdf}) has the following representation
\begin{equation}
  \label{eq:gs}
  \nu(x)=\mu([0,x]^\comp)=\int_\SS \max_{1\leq i\leq d}
  \left(\frac{a_i}{x_i}\right)\sigma(da)\,,\quad
  x\in[0,\infty]^d\setminus\{0\}\,, 
\end{equation}
where 
\begin{displaymath}
  [0,x]=\times_{i=1}^d [0,x_i]\,,\quad x=(x_1,\dots,x_d)\,,
\end{displaymath}
$\SS=\{x\in\EE:\;\|x\|=1\}$ is a sphere in $\EE=[0,\infty)^d$ with
respect to any chosen norm (from now on called the \emph{reference
  sphere} and the \emph{reference norm}) and $\sigma$ is a finite
measure on $\SS$ (called the \emph{spectral measure} of $\xi$) such
that
\begin{equation}
  \label{eq:1dim-norm}
  \int_{\SS} a_i\sigma(da)=1\,,\quad i=1,\dots,d\,.
\end{equation}
A similar representation is described in
\cite[Th.~4.3.1]{fal:hus:rei04} for the special case of $\SS$ being
the unit simplex.  

We now aim to relate the function $\nu(x)$ from (\ref{eq:gs}) to the
support function of a certain compact convex set. Recall that the
\emph{support function} of a set $M\subset\Rd$ is defined as
\begin{displaymath}
  h(M,x)=\sup\{\langle z,x\rangle:\; z\in M\}\,,
\end{displaymath}
where $\langle z,x\rangle$ is the scalar product in $\Rd$.  Let
$e_1,\dots,e_d$ be the standard orthonormal basis in $\Rd$. For every
$a=(a_1,\dots,a_d)\in\Rd$ consider the \emph{cross-polytope}
\begin{displaymath}
  \simplex{a}=\conv(\{0,a_1e_1,\dots,a_de_d\})\,,
\end{displaymath}
where $\conv(\cdot)$ denotes the convex hull of the corresponding set.
Note that $\conv(\{a_1e_1,\dots,a_de_d\})$ is a simplex.  Then
\begin{displaymath}
  h(\simplex{a},x)=h(\simplex{x},a)= \max_{1\leq i\leq d}\; (a_ix_i)
\end{displaymath}
for every $a\in\SS$ and $x\in \EE$.  For $x=(x_1,\dots,x_d)\in\EE$
write $x^*=(x_1^{-1},\dots,x_d^{-1})$. Then (\ref{eq:gs}) can be
expressed as
\begin{equation}
  \label{eq:gxk}
  \nu(x^*)=\int_\SS h(\simplex{a},x)\sigma(da)\,, \quad x\in\EE\,.
\end{equation}
The function $l(x)=\nu(x^*)$ is called the \emph{stable tail
  dependence function}, see \cite[p.~257]{beir:goeg:seg:04}.

It is well known that the arithmetic sum of support functions of two
convex compact sets $K$ and $L$ is the support function of their
Minkowski sum 
\begin{displaymath}
  K+L=\{x+y: \;x\in K, \,y\in L\}\,,
\end{displaymath}
i.e. $h(K+L,x)$ equals $h(K,x)+h(L,x)$.  Extending this idea to
integrals of support functions leads to the expectation concept for
random convex compact sets, see \cite{art:vi75} and
\cite[Sec.~2.1]{mo1}. If $X$ is a random compact convex set \cite{mo1}
such that $\|X\|=\sup\{\|x\|:\; x\in X\}$ is integrable, then the
\emph{selection expectation} (also called the Aumann expectation) of
$X$ is the set of expectations of $\E\xi$ for all random vectors $\xi$
such that $\xi\in X$ a.s.  If the underlying probability space is
non-atomic, or $X$ is a.s.  convex, then $\E X$ is the unique compact
convex set that satisfies
\begin{displaymath}
  \E h(X,x)=h(\E X,x)
\end{displaymath}
for all $x$, see \cite[Th.~II.1.22]{mo1}.

Let $\sigma_1$ be the spectral measure $\sigma$ normalised to have the
total mass 1. If $\eta$ is distributed on $\SS$ according to
$\sigma_1$, then $\simplex{\eta}$ is a random convex compact set whose
selection expectation satisfies
\begin{equation}
  \label{eq:int-sigma}
  h(\E \simplex{\eta},x)
  =\frac{1}{\sigma(\SS)}\int_\SS h(\simplex{a},x)\sigma(da) \,.
\end{equation}
Condition (\ref{eq:1dim-norm}) further implies that
\begin{equation}
  \label{eq:mar-norm}
  \sigma(\SS) h(\E \simplex{\eta},e_i)=1\,, \quad i=1,\dots,d\,.
\end{equation}
Since $h(\E \simplex{\eta},e_i)=\E h(\simplex{\eta},e_i)=\E\eta_i$, we
have $\sigma(\SS)\E\eta_i=1$ for $i=1,\dots,d$.  Together with
(\ref{eq:gs}) and (\ref{eq:cdf}) these reasons lead to the following
result.

\begin{theorem}
  \label{thr:1}
  A random vector $\xi$ is max-stable with unit Fr\'echet
  marginals if and only if its cumulative distribution function
  $F(x)=\Prob{\xi\leq x}$ satisfies
  \begin{displaymath}
    F(x) = \exp\{- c h(\E\simplex{\eta},x^*)\}\,, \quad x\in\EE\,,
  \end{displaymath}
  for a constant $c>0$ and a random vector $\eta\in\SS$ such that
  $c\E\eta=(1,\dots,1)$.
\end{theorem}

If $K=c\E\simplex{\eta}$, then
\begin{equation}
  \label{eq:fx=e-hk-x}
  F(x)=e^{-h(K,x^*)}\,, \quad x\in\EE\,.
\end{equation}
Furthermore, note that $K=\E\simplex{c\eta}$ with
$\E(c\eta)=(1,\dots,1)$. 

\begin{definition}
  \label{def:ds}
  The set $K=c\E\simplex{\eta}$ where $c>0$ and $\eta$ is a random
  vector on $\SS$ is said to be a \emph{max-zonoid}. If $\sigma_1$ is
  the distribution of $\eta$, then $\sigma=c\sigma_1$ is the
  \emph{spectral measure} of $K$. If $c\E\eta=(1,\dots,1)$, then the
  max-zonoid $K$ is called the \emph{dependency set} associated with
  the spectral measure $\sigma$ (or associated with the corresponding
  simple max-stable random vector).
\end{definition}

\begin{prop}
  \label{prop:max-z}
  A convex set $K$ is a max-zonoid if and only if there exists a
  semi-simple max-stable vector $\xi$ with cumulative distribution
  function $F(x)=e^{-h(K,x^*)}$ for all $x\in\EE$.
\end{prop}
\begin{proof}
  \textsl{Sufficiency.} A semi-simple max-stable $\xi$ can be obtained
  as $\xi=a\xi'=(a_1\xi'_1,\dots,a_d\xi'_d)$ for simple max-stable
  vector $\xi'$ and $a=(a_1,\dots,a_d)\in(0,\infty)^d$.  Let $K'$ be
  the dependency set of $\xi'$. By Theorem~\ref{thr:1},
  \begin{displaymath}
    \Prob{\xi\leq x}=\Prob{a\xi'\leq x}
    =e^{-h(K',ax^*)}=e^{-h(K,x^*)}\,,\quad x\in\EE\,,
  \end{displaymath}
  for $K=a K'=\{(a_1x_1,\dots,a_dx_d):\; (x_1,\dots,x_d)\in K\}$.

  \textsl{Necessity.} If $K$ is a max-zonoid, then $K'=a K$ is a
  dependency set for some $a\in(0,\infty)^d$. If $\xi'$ is max-stable
  with dependency set $K'$, then it is easily seen that $a\xi'$ has
  the cumulative distribution function $e^{-h(K,x^*)}$.
\end{proof}

Proposition~\ref{prop:max-z} means that each max-zonoid can be
rescaled to become a dependency set.

\begin{prop}
  \label{prop:bnd}
  A max-zonoid $K$ always satisfies 
  \begin{equation}
    \label{eq:bnd}
    \simplex{z}\subset K\subset [0,z]
  \end{equation}
  for some $z\in\EE$.
\end{prop}
\begin{proof}
  The result follows from the following bound on the support function
  of $\E\simplex{\eta}$
  \begin{displaymath}
    h(\simplex{y},x)=
    \max_{1\leq i\leq d} \E (\eta_ix_i) \leq 
    \E h(\simplex{\eta},x)\leq \E\sum_{i=1}^d \eta_ix_i=h([0,y],x)
  \end{displaymath}
  where $y=\E\eta$, so that (\ref{eq:bnd}) holds with $z=cy$.
\end{proof}

The normalisation condition (\ref{eq:mar-norm}) and (\ref{eq:bnd})
imply that the dependency set of a simple max-stable distribution
satisfies
\begin{equation}
  \label{eq:ksub}
  \simplex{(1,\dots,1)}=\conv\{0,e_1,\dots,e_d\}
  \subset K\subset [0,1]^d\,,
\end{equation}
where $\simplex{(1,\dots,1)}$ is called the unit cross-polytope.

The selection expectation of $\simplex{\eta}$ has the support function
given by 
\begin{equation}
  \label{eq:hesimpl-x=int_ss-a_1}
  h(\E\simplex{\eta},x)=\int_{\SS} \|(a_1x_1,\dots,a_dx_d)\|_\infty
  \sigma(da)\,, 
\end{equation}
where $\|\cdot\|_\infty$ is the $\ell_\infty$-norm in $\R^d$. If the
$\ell_\infty$-norm in (\ref{eq:hesimpl-x=int_ss-a_1}) is replaced by
the $\ell_1$-norm, i.e. the absolute value of the sum of the
coordinates and integration is carried over the whole sphere, then
(\ref{eq:hesimpl-x=int_ss-a_1}) yields the support function of a
\emph{zonoid}, see \cite[Sec.~3.5]{schn}. This provides one of the
reasons for calling $\E\simplex{\eta}$ a max-zonoid. Note that
max-zonoids form a sub-family of sets called $d$-zonoids in
\cite{ric82}.

It is possible to define a max-zonoid as the selection expectation of
$\simplex{\zeta}$, where $\zeta$ is any random vector in $\EE$ (not
necessarily on $\SS$). The corresponding spectral measure $\sigma$ on
$\SS$ can be found from
\begin{equation}
  \label{eq:int-fasigm-fzet}
  \int_{\SS} g(a)\sigma(da)=\E\Big[\|\zeta\| g(\frac{\zeta}{\|\zeta\|})\Big]
\end{equation}
for all integrable functions $g$ on $\SS$. Indeed,
\begin{displaymath}
  \int_{\SS} h(\simplex{u},x)\sigma(du)=
  \E[\|\zeta\| h(\simplex{\zeta/\|\zeta\|},x)]
  =\E h(\simplex{\zeta},x)\,.
\end{displaymath}
If all coordinates of $\zeta$ have the unit mean, then the selection
expectation of $\simplex{\zeta}$ becomes a dependency set.

An alternative representation of max-stable laws
\cite[Prop.~5.11]{res87} yields that
\begin{equation}
  \label{eq:mfs}
  F(x)=\exp\left\{-\int_0^1\max\left(\frac{f_1(s)}{x_1},\dots,
    \frac{f_d(s)}{x_d}\right)ds\right\}
\end{equation}
for non-negative integrable functions $f_1,\dots,f_d$ satisfying 
\begin{displaymath}
  \int_0^1 f_i(s)ds=1\,,\quad i=1,\dots,d\,.
\end{displaymath}
Thus
\begin{displaymath}
  h(K,x)=\int_0^1\max(f_1(s)x_1,\dots,f_d(s)x_d)\,ds\,,
\end{displaymath}
i.e. the dependency set $K$ is given by the selection expectation of
the cross-polytope $\simplex{f(\eta)}$, where
$f(\eta)=(f_1(\eta),\dots,f_d(\eta))$ and $\eta$ is uniformly
distributed on $[0,1]$. The corresponding spectral measure can be
found from (\ref{eq:int-fasigm-fzet}) for $\zeta=f(\eta)$.

\begin{theorem}
  \label{thr:2}
  If $d=2$, then each convex set $K$ satisfying (\ref{eq:ksub}) is the
  dependency set of a simple max-stable distribution.  If $d\geq3$,
  then only those $K$ that satisfy (\ref{eq:ksub}) and are max-zonoids
  correspond to simple max-stable distributions.
\end{theorem}
\begin{proof}
  Consider a planar convex polygon $K$ satisfying (\ref{eq:ksub}), so
  that its vertices are $a^0=e_1,a^1,\dots,a^m=e_2$ in the
  anticlockwise order.  Then $K$ equals the sum of triangles with
  vertices $(0,0),(a^{i-1}_1-a^i_1,0),(0,a^i_2-a^{i-1}_2)$ for
  $i=1,\dots,m$, where $a^i=(a^i_1,a^i_2)$. Thus (\ref{eq:gxk}) holds
  with $\sigma$ having atoms at $u_i/\|u_i\|$ with mass $\|u_i\|$
  where $u_i=(a^{i-1}_1-a^i_1,a^i_2-a^{i-1}_2)$ for $i=1,\dots,m$.
  The approximation by polytopes yields that a general convex $K$
  satisfying (\ref{eq:ksub}) can be represented as the expectation of
  a random cross-polytope and so corresponds to a simple max-stable
  distribution.

  Theorem~\ref{thr:1} implies that all max-zonoids satisfying
  (\ref{eq:ksub}) correspond to simple max-stable distributions.  It
  remains to show that not every convex set $K$ satisfying
  (\ref{eq:ksub}) is a dependency set in dimension $d\geq3$. For
  instance, consider set $L$ in $\R^3$ which is the convex hull of
  $0,e_1,e_2,e_3$ and $(2/3,2/3,2/3)$. All its 2-dimensional faces are
  triangles, so that this set is indecomposable by
  \cite[Th.~15.3]{gruen67}. Since $L$ is a polytope, but not a
  cross-polytope, it cannot be represented as a sum of
  cross-polytopes and so is not a max-zonoid. 
\end{proof}

The support function of the dependency set $K$ equals the tail
dependence function (\ref{eq:gxk}). If an estimate $\hat{l}(\cdot)$ of
the tail dependence function is given for a finite set of directions
$u_1,\dots,u_m$, it is possible to estimate $K$, e.g. as the
intersection of half-spaces $\{x\in\EE:\; \langle x,a_i\rangle \leq
\hat{l}(a_i)\}$. However, this estimate should be use very cautiously,
since the obtained polytope $K$ is not necessarily a max-zonoid in
dimensions three and more. While this approach is justified in the
bivariate case (see also \cite{hal:taj04}), in general, it is better
to use an estimate $\hat\sigma$ of the spectral measure $\sigma$ in
order to come up with an estimator of $K$ as
\begin{displaymath}
  h(\hat{K},x)=\int_{\SS} h(\simplex{a},x)\hat\sigma(da)\,.
\end{displaymath}
Being the expectation of a cross-polytope, the obtained set is
necessarily a max-zonoid.

\bigskip

The set
\begin{displaymath}
  K^o=\{x\in\EE:\; h(K,x)\leq 1\}
\end{displaymath}
is called the \emph{polar} (or \emph{dual}) set to $K$ in $\EE$, see
\cite[Sec.~1.6]{schn} for the conventional definition where $\EE$ is
replaced by $\R^d$. If $K$ is convex and satisfies (\ref{eq:ksub}),
then its polar $K^o$ is also convex and satisfies the same condition. 
Furthermore, 
\begin{align*}
  \{x\in\EE:\; F(x)\geq \alpha\}&=\{x\in\EE:\; e^{-h(K,x^*)}\geq \alpha\}\\
  &=\{x^*:\; x\in\EE,\; h(K,x)\leq -\log\alpha\}\\
  &=(-\log\alpha)\{x^*:\; x\in K^o\}\,,
\end{align*}
i.e. multivariate quantiles of the cumulative distribution function of
a simple max-stable random vector are inverted rescaled variants of
the polar set to the dependency set $K$. The level sets of
multivariate extreme values distributions have been studied in
\cite{haan:ron98}. Note that the dimension effect described in
Theorem~\ref{thr:2} restricts the family of sets that might appear as
multivariate quantiles in dimensions $d\geq3$.

\bigskip

The ordering of dependency sets by inclusion corresponds to the
stochastic ordering of simple max-stable random vectors, i.e. if
$\xi'$ and $\xi''$ have dependency sets $K'$ and $K''$ with $K'\subset
K''$, then $\Prob{\xi'\leq x}\geq\Prob{\xi''\leq x}$ for all
$x\in\EE$. 

A metric on the family of dependency sets may be used to measure the
distance between random vectors $\xi'$ and $\xi'$ with simple
max-stable distributions.  Such distance can be defined as the
Hausdorff distance between the dependency sets of $\xi$ and $\xi'$ or
any other metric for convex sets (e.g. the Lebesgue measure of the
symmetric difference or the $L_p$-distance between the support
functions).  In the spirit of the Banach-Mazur metric for convex sets
(or linear spaces), a distance between two dependency sets $K'$ and
$K''$ can be defined as
\begin{displaymath}
  m(K',K'')=\log\inf\{\prod_{i=1}^d \lambda_i:\; K'\subset\lambda K'',\;
  K''\subset\lambda K', \lambda\in(0,\infty)^d\}\,,
\end{displaymath}
where $\lambda K=\{(\lambda_1x_1,\dots,\lambda_dx_d):\;
(x_1,\dots,x_d)\in K\}$ with $\lambda=(\lambda_1,\dots,\lambda_d)$.
If $\xi'$ and $\xi''$ have dependency sets $K'$ and $K''$
respectively, then $m(K',K'')$ is the logarithm of the smallest value
of $(\lambda_1\cdots\lambda_d)$ such that $\xi'$ is stochastically
smaller than $\lambda\xi''$ and $\xi''$ is stochastically smaller than
$\lambda\xi'$. For instance, the distance between the unit
cross-polytope and the unit square (for $d=2$) is $\log 4$, which is
the largest possible distance between two simple bivariate max-stable
laws.

\section{Norms associated with max-stable distributions}
\label{sec:norm-associated-with}

Note that the support function of a compact set $L$ is sublinear, i.e.
it is homogeneous and subadditive. If $L$ is convex symmetric and
contains the origin in its interior, then its support function
$h(L,x)$ defines a norm in $\Rd$.  Conversely, every norm defines a
symmetric convex compact set in $\Rd$ with the origin in its
interior, see \cite[Th.~15.2]{roc70}.

Let $K$ be a convex set satisfying (\ref{eq:ksub}). The corresponding
norm $\|\cdot\|_K$ can be defined as the support function of the set
$L$ obtained as the union of all symmetries of $K$ with respect to
coordinate planes, i.e.
\begin{displaymath}
  \|x\|_K=h(L,x)=h(K,|x|)\,, \quad x\in\R^d\,,
\end{displaymath}
where $|x|=(|x_1|,\dots,|x_d|)$. The norm $\|x\|_K$ is said to be
generated by the max-zonoid $K$. Note that the origin belongs to the
interior of $L$ and $\|x\|_K=h(K,x)$ for $x\in\EE$.  The following
result shows that distributions of max-stable vectors correspond to
norms generated by max-zonoids.

\begin{theorem}
  \label{thr:norm}
  Let $\|\cdot\|$ be a norm on $\R^d$. The function
  \begin{equation}
    \label{eq:fx=e-x-quad}
    F(x)=\exp\{-\|x^*\|\}\,,\quad x\in\EE\,,
  \end{equation}
  is the cumulative distribution function of a random vector $\xi$ in
  $\EE$ if and only if $\|x\|=h(K,|x|)$ is the norm generated by a
  max-zonoid $K$. In this case the random vector $\xi$ is necessarily
  semi-simple max-stable.
\end{theorem}
\begin{proof}
  \textsl{Sufficiency.} If $K$ is a max-zonoid,
  Proposition~\ref{prop:max-z} implies that there exists a semi-simple
  max-stable vector $\xi$ with cumulative distribution function
  $e^{-h(K,x^*)}=e^{-\|x^*\|}$.

  \textsl{Necessity.} If (\ref{eq:fx=e-x-quad}) is the cumulative
  distribution function of a random vector $\xi$, then
  \begin{displaymath}
    \Prob{\xi^{(1)}\vee\cdots\vee\xi^{(n)}\leq x}
    =e^{-n\|x^*\|}=e^{-\|(n^{-1}x)^*\|}
    =\Prob{\xi\leq n^{-1}x}
  \end{displaymath}
  for all $x\in\EE$ and i.i.d. copies $\xi^{(1)},\dots,\xi^{(n)}$ of
  $\xi$. Thus, $\xi$ is necessarily semi-simple max-stable.
  Proposition~\ref{prop:max-z} implies that (\ref{eq:fx=e-x-quad})
  holds with the norm generated by a max-zonoid $K$.
\end{proof}

The space $\R^d$ with the norm $\|\cdot\|_K$ becomes a
finite-dimensional normed linear space, also called the
\emph{Minkowski space}, see \cite{thom96}.  If this is an inner
product space, then the norm is necessarily Euclidean. Indeed, if $K$
is the intersection of a centred ellipsoid with $\EE$ and satisfies
(\ref{eq:ksub}), then this ellipsoid is necessarily the unit ball.

\medskip

Another common way to standardise the marginals of a multivariate
extreme value distribution is to bring them to the reverse exponential
distribution (or unit Weibull distribution), see
\cite[Sec.~4.1]{fal:hus:rei04}. In this case, the cumulative
distribution function turns out to be
\begin{displaymath}
  F(x)=e^{-\|x\|_K}\,\quad x\in(-\infty,0]^d\,.
\end{displaymath}
The fact that every max-stable distribution with reverse exponential
marginals gives rise to a norm has been noticed in
\cite[p.~127]{fal:hus:rei04}, however without giving a
characterisation of these norms.

\bigskip

Note that the norm of $\|x\|_K$ can be expressed as
\begin{displaymath}
  \|x\|_K=\|x\|\; \|u_x\|_K\,,
\end{displaymath}
where $u_x=x/\|x\|$ belongs to the reference sphere $\SS$. If the
reference norm is $\ell_1$, then $\SS$ is the unit simplex and the
norm $\|u_x\|_K$ of $u=(t_1,\dots,t_{d-1},1-t_1-\cdots-t_{d-1})\in\SS$
can be represented as a function $A(t_1,\dots,t_{d-1})$ of
$t_1,\dots,t_{d-1}\geq0$ such that $t_1+\cdots+t_{d-1}\leq1$. If
$d=2$, then $A(t)$, $0\leq t\leq1$, is called the \emph{Pickands
  function}, see \cite{kot:nad00} and for the multivariate case also
\cite{fal:reis05,kot:nad00}. In general, the norm $\|u\|$, $u\in\SS$,
is an analogue of the Pickands function.

\begin{example}
  \label{ex:basic}
  The dependency set $K$ being the unit cube $[0,1]^d$ (so that
  $\|x\|_K$ is the $\ell_1$-norm) corresponds to the independence
  case, i.e.  independent coordinates of $\xi=(\xi_1,\dots,\xi_d)$.
  The corresponding spectral measure allocates unit atoms to the
  points from the coordinate axes.
  
  Furthermore, $K$ being the unit cross-polytope (so that $\|x\|_K$ is
  the $\ell_\infty$-norm) gives rise to the random vector
  $\xi=(\xi_1,\dots,\xi_1)$ with all identical $\Phi_1$-distributed
  coordinates, i.e. the completely dependent random vector. The
  corresponding spectral measure has its only atom at the point from
  $\SS$ having all equal coordinates.  Note that the unit cube is dual
  (or polar) set to the unit cross-polytope.
\end{example}

\begin{example}
  \label{ex:ball-ellipse}
  The $\ell_p$-norm $\|x\|_p$ with $p\geq1$ generates the symmetric
  logistic distribution \cite[(9.11)]{beir:goeg:seg:04} with parameter
  $\alpha=1/p$. The strength of dependency increases with $p$.
\end{example}

\begin{example}
  \label{ex:gauss}
  A useful family of simple max-stable bivariate distribution appears
  if the functions $f_1,f_2$ in (\ref{eq:mfs}) are chosen to be the
  density functions of normal distributions, see
  \cite[Sec.~3.4.5]{kot:nad00} and \cite[p.~309]{beir:goeg:seg:04}. It
  is shown in \cite{hues:reis89} that these distributions appear as
  limiting distributions for maxima of bivariate i.i.d. Gaussian
  random vectors.  The corresponding norm (which we call the
  H\"usler-Reiss norm) is given by
  \begin{displaymath}
    \|x\|_K= x_1\Phi(\lambda+\frac{1}{2\lambda}\log \frac{x_1}{x_2})
    +x_2\Phi(\lambda-\frac{1}{2\lambda}\log \frac{x_1}{x_2})\,,
  \end{displaymath}
  where $\lambda\in[0,\infty]$. The cases $\lambda=0$ and
  $\lambda=\infty$ correspond to complete dependence and independence,
  respectively.  
\end{example}

\begin{example}
  \label{ex:neg-lo}
  The bivariate symmetric negative logistic distribution
  \cite[p.~307]{beir:goeg:seg:04} corresponds to the norm given by
  \begin{displaymath}
   \|x\|_K=\|x\|_1-\lambda\|x\|_p\,,
  \end{displaymath}
  where $\lambda\in[0,1]$ and $p\in[-\infty,0]$.
\end{example}

\section{Spectral measures}
\label{sec:expon-meas-spectr}

Since the dependency set determines uniquely the distribution of a
simple max-stable random vector, there is one to one correspondence
between dependency sets and normalised spectral measures. It is
possible to extend this correspondence to max-zonoids on one side and
all finite measures on $\SS$ on the other one, since both uniquely
identify semi-simple max-stable distributions. While the spectral
measure depends on the choice of the reference norm, the dependency
set remains the same whatever the reference norm is.

It is shown in \cite{col:taw91} that the densities of the spectral
measure on the reference simplex $\SS=\conv\{e_1,\dots,e_d\}$ can be
obtained by differentiating the function $\nu(x)=\mu([0,x]^\comp)$ for
the exponent measure $\mu$.  This is possible if the spectral measure
is absolutely continuous with respect to the surface area measures on
relative interiors of all faces of the simplex and, possibly, has
atoms at the vertices of $\SS$.  Following the proof of this fact in
\cite[Sec.~8.6.1]{beir:goeg:seg:04}, we see that
\begin{displaymath}
  \lim_{z_j\to0,\; j\notin A} D_A \nu(z)
  =(-1)^{|A|-1}D_A\mu(\{x\in\EE:\; x_j>z_j, j\in A;\, x_j=0, j\notin
  A\})\,, 
\end{displaymath}
where $A\subset\{1,\dots,d\}$, $|A|$ is the cardinality of $A$, and
$D_A$ denotes the mixed partial derivative with respect to the
coordinates with numbers from $A$.

The derivatives of $\nu$ can be expressed by means of the derivatives
of the stable tail dependence function $l(z)=\nu(z^*)=\|z\|_K$. Indeed,
\begin{displaymath}
  D_A \nu(z)=D_A l(z^*)(-1)^{|A|}\prod_{j\in A} z_j^{-2}\,.
\end{displaymath}
Thus, the densities of the exponent measure can be found from 
\begin{multline*}
  D_A\mu(\{x\in\EE:\; x_j\leq z_j, j\in A;\, x_j=0, j\notin
  A\})= (-1)^{|A|-1}\lim_{z_j\to0,\; j\notin A} 
  D_A l(z^*)\prod_{j\in A} z_j^{-2}\,.
\end{multline*}
In particular, the density of $\mu$ in the interior of $\EE$ can be
found from the $d$th mixed partial derivative of the norm as
\begin{displaymath}
  f(z)=(-1)^{d-1} \frac{\partial^d l}{\partial z_1\cdots\partial
    z_d}(z^*) \prod_{i=1}^d z_i^{-2}\,, \quad z\in(0,\infty)^d\,.
\end{displaymath}
After decomposing these densities into the radial and directional
parts, it is possible to obtain the spectral measure by
\begin{displaymath}
  \sigma(G)=\mu(\{tu:\; u\in G,\, t\geq1\})
  =\int_{\{tu:\; u\in G,\, t\geq1\}} f(z)dz
\end{displaymath}
for every measurable $G$ from the relative interior of $\SS$. The
relationship between spectral measures on two different reference
spheres is given in \cite[p.~264]{beir:goeg:seg:04}.

\begin{prop}
  \label{prop:deriv}
  A $d$-times continuously differentiable function $l(x)$, $x\in E$,
  is the tail dependency function of a simple max-stable distribution
  if and only if $l$ is sublinear, takes value 1 on all basis vectors,
  and all its mixed derivatives of even orders are non-positive and of
  odd orders are non-negative.
\end{prop}
\begin{proof}
  The necessity follows from Theorem~\ref{thr:1} and the
  non-negativity condition on the exponent measure $\mu$. In the other
  direction, the sublinearity property implies that $l$ is the support
  function of a certain convex set $K$, see \cite[Th.~1.7.1]{schn}.
  The condition on the sign of mixed derivatives yields that the
  corresponding densities of $\mu$ are non-negative, i.e. $K$ is the
  max-zonoid corresponding to a certain spectral measure.
\end{proof}

In the planar case, \cite[Th.~1.7.2]{schn} implies that the second
mixed derivative of the (smooth) support function is always
non-positive. Accordingly, all smooth planar convex sets satisfying
(\ref{eq:ksub}) are dependency sets.

A number of interesting measures on the unit sphere appear as
curvature measures of convex sets \cite[Sec.~4.2]{schn}. A complete
interpretation of these curvature measures is possible in the planar
case, where the curvature measure becomes the length measure. The
\emph{length measure} $S_1(L,A)$ generated by a smooth set $L$
associates with every measurable $A\subset\Sphere[1]$ the
1-dimensional Hausdorff measure of the boundary of $L$ with unit
normals from $A$. The length measure for a general $L$ is defined by
approximation. Recall that $\Sphere[1]_+$ is the part of the unit
circle lying in the first quadrant.

\begin{theorem}
  \label{thr:2d}
  A measure $\sigma$ on $\Sphere[1]_+$ is the spectral measure of a
  simple max-stable random vector $\xi$ with dependency set $K$ if and
  only if $\sigma$ is the restriction on $\Sphere[1]_+$ of the length
  measure generated by $\check{K}=\{(x_1,x_2):\; (x_2,x_1)\in K\}$
  with $K$ satisfying (\ref{eq:ksub}).  
\end{theorem}
\begin{proof}
  \textsl{Sufficiency.}
  Consider a planar convex set $K$ satisfying (\ref{eq:ksub}). Let
  $\sigma(da)$ be the length measure of $L=\check{K}$, i.e.  the
  first-order curvature measure $S_1(L,da)$.  Then
  \begin{displaymath}
    \int_{\Sphere[1]} h(\simplex{a},x)\sigma(da)
    =\int_{\Sphere[1]} h(\simplex{x},a)S_1(L,da)
    =2 V(\simplex{x},L)\,,
  \end{displaymath}
  where $V(\simplex{x},L)$ denotes the mixed volume (the mixed area in
  the planar case) of the sets $\simplex{x}$ and $L$, i.e. 
  \begin{displaymath}
    2V(\simplex{x},L)=\mes_2(L+\simplex{x})-\mes_2(L)-\mes_2(\simplex{x})\,,
  \end{displaymath}
  see \cite[Sec.~5.1]{schn}. Because of (\ref{eq:ksub}), the integral
  over the full circle $\Sphere[1]$ with respect to $\sigma$ coincides
  with the integral over $\Sphere[1]\cap[0,\infty)^2=\Sphere[1]_+$.
  It remains to show that $2V(\simplex{x},L)$ equals $h(K,x)$.
  If $z=(z_1,z_2)\in L$ is any support point of $L$ in direction
  $x=(x_1,x_2)$, i.e.  $h(L,x)=\langle z,x\rangle$, then
  \begin{displaymath}
    2V(\simplex{x},L)=z_1x_2+z_2x_1=h(K,x)\,.
  \end{displaymath}

  An alternative proof follows the construction from
  Theorem~\ref{thr:2}. Indeed, a polygonal $K$ can be obtained as the
  sum of triangles. A triangle $\simplex{(t,s)}$ with vertices
  $(0,0),(t,0)$ and $(0,s)$ corresponds to the spectral measure having
  the atom at $(t,s)c^{-1}$ with mass $c=\sqrt{t^2+s^2}$ and therefore
  coincides with the length measure of
  $\simplex{(s,t)}=\check{\Delta}_{(t,s)}$ restricted onto
  $\Sphere[1]_+$. Since the spectral measure of $K$ is the sum of
  spectral measures of these triangles, it can be alternatively
  represented as the sum of the length measures. A general $K$ can be
  then approximated by polygons. 
  
  \textsl{Necessity.} Assume that a measure $\sigma$ on
  $\Sphere[1]_+$ is the spectral measure of a simple max-stable law
  with dependency set $K$. If now $\sigma'$ is chosen to be the length
  measure of $\check{K}$ restricted onto $\Sphere[1]_+$, then
  $\sigma'$ generates the max-zonoid $K$. Finally, $\sigma=\sigma'$ by
  the uniqueness of the spectral measure.
\end{proof}

Theorem~\ref{thr:2d} implies that the length of the boundary of $K$
inside $(0,\infty)^2$ equals the total mass of the spectral measure on
$\Sphere[1]_+$.  Given (\ref{eq:ksub}), an obvious bound on this
boundary length implies that this total mass lies between $\sqrt{2}$
and $2$.

The total mass of the spectral measure on the reference \emph{simplex}
has a simple geometric interpretation.  Assume that the reference norm
is $\ell_1$, i.e.  $\|x\|=x_1+\cdots+x_d$ for $x\in\EE$. If
$\eta\in\SS$, then $\E\eta_1+\cdots+\E\eta_d=1$, so that the
$\ell_1$-norm of $\E\eta$ is $1$. Since $c\E\eta=(1,\dots,1)$ in
Theorem~\ref{thr:1}, we have $c=d$, i.e. the spectral measure has the
total mass $d$.

The weak convergence of simple max-stable random vectors can be
interpreted as convergence of the corresponding max-zonoids.

\begin{theorem}
  \label{prop:wc}
  Let $\xi,\xi_1,\xi_2,\ldots$ be a sequence of simple max-stable
  random vectors with spectral measures
  $\sigma,\sigma_1,\sigma_2,\ldots$ and dependency sets
  $K,K_1,K_2,\ldots$. Then the following statements are equivalent.
  \begin{description}
  \item[(i)] $\xi_n$ converges in distribution to $\xi$;
  \item[(ii)] $\sigma_n$ converges weakly to $\sigma$;
  \item[(iii)] $K_n$ converges in the Hausdorff metric to $K$.
  \end{description}
\end{theorem}
\begin{proof}
  The equivalence of (i) and (ii) is well known, see
  \cite[Cor.~6.1.15]{haan:fer06}. 

  The weak convergence of $\sigma_n$, the continuity of
  $h(\simplex{a},x)$ for $a\in\SS$ and (\ref{eq:int-sigma}) imply that
  the support function of $K_n$ converges pointwisely to the support
  function of $K$.  Because dependency sets are contained inside the
  unit cube and so are uniformly bounded, their convergence in the
  Hausdorff metric is equivalent to the pointwise convergence of their
  support functions.

  The Hausdorff convergence of $K_n$ to $K$ implies the pointwise
  convergence of their support functions and so the pointwise
  convergence of the cumulative distribution functions given by
  (\ref{eq:fx=e-hk-x}). The latter entails that $\xi_n$ converges in
  distribution to $\xi$.
\end{proof}

A random vector $\zeta\in\EE$ belongs to the \emph{domain of
  attraction} of a simple max-stable distribution if and only if the
measure
\begin{equation}
  \label{eq:sigma-s}
  \sigma_s(A)=s\Prob{\frac{\zeta}{\|\zeta\|}\in A\,,\,
    \|\zeta\|\geq s}\,, \quad A\subset\SS\,,
\end{equation}
converges weakly as $s\ti$ to a finite measure on $\SS$, which then
becomes the spectral measure of the limiting random vector, see
\cite[(8.95)]{beir:goeg:seg:04}. The equivalence of (ii) and (iii) in
Theorem~\ref{prop:wc} implies the following result.

\begin{prop}
  \label{prop:weak-conv}
  A random vector $\zeta\in\EE$ belongs to the domain of attraction of
  a simple max-stable random vector $\xi$ with spectral measure
  $\sigma$ if and only if the max-zonoids generated by $\sigma_s$ from
  (\ref{eq:sigma-s}) converge in the Hausdorff metric as $s\ti$ to the
  max-zonoid generated by $\sigma$.
\end{prop}

\section{Copulas and association}
\label{sec:copulas-association}

The dependency structure of a distribution with fixed marginals can be
explored using the copula function $C$ defined on $\II^d=[0,1]^d$ by
the following equation
\begin{displaymath}
  F(x)=F(x_1,\dots,x_d)=C(F_1(x_1),\dots,F_d(x_d))\,,
\end{displaymath}
where $F_1,\dots,F_d$ are the marginals of $F$, see \cite{nel06}. In
case of a simple max-stable distribution, we obtain
\begin{equation}
  \label{eq:cop}
  C(u_1,\dots,u_d)=\exp\{-\|(-\log u_1,\dots,-\log u_d)\|_K\}\,.
\end{equation}

\begin{theorem}
  \label{th:cop}
  The function (\ref{eq:cop}) is a copula function if and only if $K$
  is a max-zonoid.
\end{theorem}
\begin{proof}
  The sufficiency is trivial, since the right-hand side of
  (\ref{eq:cop}) can be used to construct a max-stable distribution.
  In the other direction, we can substitute into $C$ the
  $\Phi_1$-distribution functions, i.e. $u_i=e^{-1/x_i}$.  This yields
  a multivariate cumulative distribution function given by
  $F(x)=e^{-\|x^*\|_K}$.  The result then follows from
  Theorem~\ref{thr:norm}.
\end{proof}

Note that (\ref{eq:cop}) in the bivariate case appears in
\cite{falk06}.  A rich family of copulas consists of the
\emph{Archimedean copulas} that in the bivariate case satisfy
$\phi(C(x_1,x_2))=\phi(x_1)+\phi(x_2)$ for a strictly decreasing
continuous function $\phi$ and all $x_1,x_2\in[0,1]$, see
\cite[Ch.~4]{nel99}.  Using (\ref{eq:cop}), it is easy to see that in
this case $\psi(\|(x_1,x_2)\|_K)=\psi(x_1)+\psi(x_2)$ for a monotone
increasing continuous function $\psi$ and all $x_1,x_2\geq0$.  It is
known (see \cite[Th.~4.5.2]{nel06} and \cite{gen:riv89}) that all
Archimedean copulas that correspond to max-stable distributions are
so-called Gumbel copulas, where $\psi(t)=t^p$. Thus, (\ref{eq:cop}) is
an Archimedean copula if and only if $K$ is $\ell_p$-ball with
$p\in[1,\infty]$, see Example~\ref{ex:ball-ellipse}.

\bigskip

The bivariate copulas are closely related to several association
concepts between random variables, see \cite{nel91}. The
\emph{Spearman correlation coefficient} is expressed as
$\rho_S=12J-3$, where
\begin{align*}
  J &=\int_0^1\int_0^1 C(u_1,u_2)du_1du_2
  =\int_0^1\int_0^1 e^{-\|(-\log u_1,-\log u_2)\|_K}du_1du_2\\
  &=\int_0^\infty\int_0^\infty
  e^{-\|(x_1,x_2)\|_K}e^{-x_1-x_2}dx_1dx_2\\
  &=\frac{1}{4}\int_(0,\infty)^2 e^{-\|x\|_L}dx\,,
\end{align*}
and
\begin{displaymath}
  L=\frac{1}{2}(K+\II^2)\,.
\end{displaymath}
It is possible to calculate $J$ by changing variables $x=r(t,1-t)$
with $r\geq0$ and $t\in[0,1]$, which leads to the following known
expression
\begin{displaymath}
  J=\frac{1}{4}\int_0^1 \frac{1}{\|(t,1-t)\|_L^2} dt 
  =\int_0^1 \frac{1}{(1+\|(t,1-t)\|_K)^2} dt\,, 
\end{displaymath}
see \cite{huer03}.  The following proposition is useful to provide
another geometric interpretation of $\rho_S$ and also an alternative
way to compute $J$.

\begin{prop}
  \label{prop:int}
  If $L$ is a convex set in $\R^d$, then 
  \begin{displaymath}
    \int_{[0,\infty)^d} e^{-h(L,x)}dx=\Gamma(d+1)\mes_d(L^o)\,,
  \end{displaymath}
  where $\mes_d(\cdot)$ is the $d$-dimensional Lebesgue measure, $L^o$
  is the polar set to $L$ and $\Gamma$ is the Gamma function.
\end{prop}
\begin{proof}
  The proof follows the argument mentioned in \cite[p.~2173]{vi96w}.
  Let $\zeta$ be the exponentially distributed random variable of mean
  $1$. Then
  \begin{align*}
    \int_{[0,\infty)^d} e^{-h(L,x)}dx
    &=\E\int_{[0,\infty)^d} \one_{\zeta\geq h(L,x)}dx\\
    &=\E\mes_d(\{x\in\EE:\; h(L,x)\leq \zeta\})\\
    &=\E \mes_d(\zeta L^o)=\mes_d(L^o) \E\zeta^d\,.
  \end{align*}
  It remains to note that $\E\zeta^d=\Gamma(d+1)$.
\end{proof}

Thus, in the planar case
\begin{displaymath}
  \rho_S=3(2\mes_2(L^o)-1)\,.
\end{displaymath}
As a multivariate extension, an affine function
$\rho_S=c(\mes_d(L^o)-a)$ of the $d$-dimensional volume of $L^o$ may
be used to define the Spearman correlation coefficient for a
$d$-dimensional max-stable random vector with unit Fr\'echet
marginals. By considering the independent case $L=K=\II^d$ with
$\rho_S=0$ and using the fact that the volume of the unit
cross-polytope $L^o$ is $(d!)^{-1}$, we see that
$\rho_S=c(d!\mes_d(L^o)-1)$ for some constant $c>0$. The choice
$c=(d+1)/(2^d-d-1)$ ensures that $\rho_S=1$ in the totally dependent
case, where $\mes_d(L^o)=2^d/(d+1)!$.

\bigskip

The \emph{Kendall correlation coefficient} of a bivariate copula $C$
is given by
\begin{align*}
  \tau&=4\int_0^1\int_0^1 C(z_1,z_2)dC(z_1,z_2) -1 \\
  &=1-4\int_0^1\int_0^1 \frac{\partial}{\partial z_1} C(z_1,z_2)
  \frac{\partial}{\partial z_2}  C(z_1,z_2) dz_1dz_2\,,
\end{align*}
see \cite[(2.3)]{nel91}.  By (\ref{eq:cop}), the partial derivatives
of $C$ can be expressed using partial derivatives of the support
function of $K$. The directional derivative of the support function
$h(K,x)$ at point $x$ in direction $u$ is given by $h(F(K,x),u)$,
where
\begin{displaymath}
  F(K,x)=\{y\in K:\; \langle y,x\rangle= h(K,x)\}
\end{displaymath}
is the \emph{support set} of $K$ in direction $x$, see
\cite[Th.~1.7.2]{schn}. Thus the partial derivatives of $h(K,x)$ are
given by
\begin{displaymath}
  \frac{\partial h(K,x)}{\partial x_i}=h(F(K,x),(1,0))=y_i(K,x)\,,
  \quad i=1,2\,,
\end{displaymath}
where $y_1(K,x)$ and $y_2(K,x)$ are respectively the maximum first and
second coordinates of the points from $F(K,x)$. If the dependency set
$K$ is strictly convex in $(0,\infty)^2$, i.e. the boundary of $K$
inside $(0,\infty)^2$ does not contain any segment, then
$F(K,x)=\{(y_1(K,x),y_2(K,x))\}$ is a singleton for all $x\in\EE$.  In
this case denote
\begin{displaymath}
  y(K,x)=y_1(K,x)y_2(K,x)\,.
\end{displaymath}
By using (\ref{eq:cop}) and changing variables we arrive at
\begin{displaymath}
  \tau=1-4\int_{[0,\infty)^2} e^{-2\|x\|_K}y(K,x)dx\,. 
\end{displaymath}
The fact that $y(K,tx)=y(K,x)$ and a similar argument to
Proposition~\ref{prop:int} yield that
\begin{equation}
  \label{eq:tau=1-2int_ko-yk}
  \tau=1-2\int_{K^o} y(K,x)dx\,.
\end{equation} 
For instance, $\tau=1/2$ if $\xi$ has the logistic distribution with
parameter $1/2$, i.e. $\|\cdot\|_K$ is the Euclidean norm.  By
changing variables $x=(t,1-t)r$, we arrive at
\begin{displaymath}
  \tau=1-\int_0^1 \frac{y_1(K,(t,1-t))y_2(K,(t,1-t))}{\|(t,1-t)\|_K^2}
  dt\,,
\end{displaymath}
which also corresponds to \cite[Th.~3.1]{huer03}.  

\bigskip

The \emph{Pearson correlation coefficient} for the components of a
bivariate simple max-stable random vector is not defined, since the
unit Fr\'echet marginals are not integrable. However it is possible to
compute it for the inverted coordinates of $\xi$.

\begin{prop}
  \label{prop:icorr}
  If $\xi$ is a simple max-stable bivariate random vector, then
  $\E(\xi_1^{-1}\xi_2^{-1})=2\mes_2(K^o)$, and the covariance between
  $\xi_1^{-1}$ and $\xi_2^{-1}$ is $2\mes_2(K^o)-1$.
\end{prop}
\begin{proof}
  Integrating by parts, it is easy to see that
  \begin{displaymath}
    \E(\xi_1^{-1}\xi_2^{-1})=
    \int_0^\infty\int_0^\infty F(x^*)dx_1dx_2
    =\int_\EE e^{-h(K,x)}dx\,.
  \end{displaymath}
  The result follows from Proposition~\ref{prop:int} and the fact that
  $\E(\xi_1^{-1})=\E(\xi_2^{-1})=1$.
\end{proof}

Proposition~\ref{prop:icorr} corresponds to the formula 
\begin{displaymath}
  \rho=\int_0^1 \frac{1}{\|(t,1-t)\|_K^2} dt -1\,.
\end{displaymath}
for the covariance obtained in \cite{taw88} for the exponential
marginals. 

Extending this concept for the higher-dimensional case, we see that
the covariance matrix of $\xi^*$ is determined by the areas of polar
sets to the 2-dimensional projections of $K$ and
\begin{displaymath}
  \rho=\frac{d!\mes_d(K^o)-1}{d!-1}
\end{displaymath}
can be used to characterise the multivariate dependency of a simple
$d$-dimensional max-stable random vector $\xi$, so that $\rho$ varies
between zero (complete independence) and $1$ (complete dependence).

\begin{example}
  \label{ex:lp-vol}
  Assume that $\|x\|_K=\|x\|_p$ is the $\ell_p$-norm with $p\geq1$,
  i.e. the corresponding $\xi$ has the logistic distribution with
  parameter $\alpha=1/p$, see Example~\ref{ex:ball-ellipse}. The
  volume of the $\ell_p$-ball $\{x\in\R^d:\; \|x\|_p\leq 1\}$ equals
  \begin{displaymath}
    v_d(p)=\frac{(2\Gamma(1+1/p))^d}{\Gamma(1+d/p)}\,,
  \end{displaymath}
  see \cite[p.~11]{pis89}. Thus, the volume of $K^o$ is $2^{-d}v_d(p)$
  and the multivariate dependency of $\xi$ can be described by
  \begin{displaymath}
    \rho=\frac{1}{d!-1}\left(d!
      \frac{(\Gamma(1+1/p))^d}{\Gamma(1+d/p)}-1\right)\,.
  \end{displaymath}
  If $d=2$, then $\rho=\alpha B(\alpha,\alpha)-1$ with $B$ being the
  Beta-function.
\end{example}

\bigskip

The \emph{tail dependency index} for $\xi=(\xi_1,\xi_2)$ with
identical marginal distributions supported by the whole positive
half-line is defined as
\begin{displaymath}
  \chi=\lim_{t\ti} \Prob{\xi_2>t|\xi_1>t}\,.
\end{displaymath}
An easy argument shows that $\chi=2-\|(1,1)\|_K$ if $\xi$ has a simple
max-stable distribution with dependency set $K$, cf
\cite{col:hef:twan99}.

\medskip

It is easy to see that $\xi$ has all independent coordinates if and
only if $\|(1,\dots,1)\|_K=d$ and the completely dependent coordinates
if and only if $\|(1,\dots,1)\|_K=1$, cf \cite{tak94} and
\cite[p.~266]{beir:goeg:seg:04}.  It is well known
\cite[p.~266]{beir:goeg:seg:04} that the pairwise independence of the
coordinates of $\xi$ implies the joint independence. Indeed, the
spectral measure of the set $u\in\SS$ such that at least two
coordinates of $u$ are positive is less than the sum of
$\sigma\{u\in\SS:\; u_i>0,u_j>0\}$ over all $i\neq j$. Each of these
summands vanishes, since
\begin{displaymath}
  x_i+x_j-\int_{\SS} \max_{1\leq k\leq d}(u_kx_k) \sigma(du)
  =\int_{\SS}((u_ix_i+u_jx_j)-(u_ix_i\vee u_jx_j))\sigma(du)=0
\end{displaymath}
by the pairwise independence, where $x$ has all vanishing coordinates
apart from $x_i$ and $x_j$. This leads to the following property of
max-zonoids. 

\begin{prop}
  \label{prop:mzp}
  If $K$ is a max-zonoid with all its two-dimensional projections
  being unit squares, then $K$ is necessarily the unit cube.
\end{prop}

\section{Complete alternation and extremal coefficients}
\label{sec:compl-altern-extr}

Consider a numerical function $f$ defined on a semigroup $S$ with a
commutative binary operation $+$. For $n\geq1$ and $x_1,\dots,x_n\in
S$ define the following successive differences
\begin{align*}
  \Delta_{x_1}f(x)&=f(x)-f(x+x_1)\,,\\
  &\cdots\\
  \Delta_{x_n}\cdots\Delta_{x_1}f(x)
  &=\Delta_{x_{n-1}}\cdots\Delta_{x_1}f(x)-
  \Delta_{x_{n-1}}\cdots\Delta_{x_1}f(x+x_n)\,.
\end{align*}
The function $f$ is said to be \emph{completely alternating} (resp.
\emph{monotone}) if all these successive difference are non-positive
(resp.  non-negative), see \cite[Sec.~4.6]{ber:c:r} and
\cite[Sec.~I.1.2]{mo1}. We will use these definitions in the following
cases: $S$ is the family of closed subsets of $\R^d$ with the union
operation, $S$ is $\R^d$ or $\EE$ with coordinatewise minimum or
coordinatewise maximum operation. Then we say shortly that the
function is max-completely alternating or monotone (resp.
min-completely or union-completely).

Every cumulative distribution function $F$ is min-completely monotone.
This is easily seen by considering the random set $X=\{\xi\}$ where
$\xi$ has the distribution $F$, then noticing that $F(x)=\Prob{X\cap
  L_x=\emptyset}=Q(L_x)$ with $L_x$ being the complement to
$x+(-\infty,0)^d$ is the avoiding functional of $X$, and finally using
the fact that 
\begin{displaymath}
  F(\min(x,y))=\Prob{X\cap (L_x\cup L_y)=\emptyset}=Q(L_x\cup L_y)
\end{displaymath}
together with the union-complete monotonicity of $Q$, see
\cite[Sec.~1.6]{mo1}.

\begin{theorem}
  \label{thr:cam}
  A convex set $K\subset\EE$ is a max-zonoid if and only if $h(K,x)$
  is a max-completely alternating function of $x$.
\end{theorem}
\begin{proof}
  It follows from \cite[Prop.~4.6.10]{ber:c:r} that a function $f(x)$
  on a general semigroup is completely alternating if and only if
  $F(x)=e^{-tf(x)}$ is completely monotone for all $t>0$.  Since
  $(\min(x,y))^*=\max(x^*,y^*)$ for $x,y\in\EE$, the function
  $x\mapsto f(x^*)$ is max-completely alternating on $\EE$ if and only
  if $x\mapsto f(x)$ is min-completely alternating.
  
  If $K$ is the max-zonoid corresponding to a simple max-stable random
  vector with distribution function $F(x)=e^{-h(K,x^*)}$, then $F^t$
  is also a cumulative distribution function (and so is min-completely
  monotone) for each $t>0$. Thus, $h(K,x^*)$ is min-completely
  alternating, whence $h(K,x)$ is max-completely alternating.
  
  In the other direction, if $h$ is max-completely alternating, then
  $F(x)=e^{-h(K,x^*)}$ is min-completely monotone, whence it is a
  cumulative distribution function. The corresponding law is
  necessarily semi-simple max-stable, so that $K$ is indeed a
  max-zonoid. 
\end{proof}

Theorem~\ref{thr:cam} can be compared with a similar characterisation
of classical zonoids, where the complete alternation of $h(K,x)$ and
monotonicity of $e^{-h(K,x)}$ are understood with respect to the
vector addition on $\R^d$, see \cite[p.~194]{schn}.

\medskip

The remainder of this section concerns extensions for the support
function defined on a finite subset of $\EE$.  

\begin{theorem}
  \label{thr:consist}
  Let $M$ be a finite set in $\EE$, which is closed with respect to
  coordinatewise maxima, i.e. $u\vee v\in M$ for all $u,v\in M$.
  Assume that for each $u,v\in M$, we have $tu\leq v$ if and only if
  $u\leq v$ and $t\leq1$.  Then a non-negative function $h$ on $M$ can
  be extended to the support function of a max-zonoid if and only if
  $h$ is max-completely alternating on $M$.
\end{theorem}
\begin{proof}
  The necessity trivially follows from Theorem~\ref{thr:cam}. To prove
  the sufficiency we explicitly construct (following the ideas of
  \cite{sch:taw02}) a max-stable random vector $\xi$ such that the
  corresponding norm coincides with the values of $h$ on the points
  from $M$.
  
  For any set $A\subset M$, let $\vee A$ denote the coordinatewise
  maximum of $A$. Furthermore, define $T(A)=h(\vee A)/h(\vee M)$.
  Since $h$ is max-completely alternating, $T$ is union-completely
  alternating on subsets of $M$.  The Choquet theorem
  \cite[Th.~I.1.13]{mo1} implies that a union-completely alternating
  function on a discrete set is the capacity functional
  $T(A)=\Prob{X\cap A\neq\emptyset}$ of a random closed set $X\subset
  M$. Define $c_u=h(\vee M)\Prob{\vee X=u}$ for $u\in M$.
  
  Let $\zeta_u$, $u\in M$, be the family of i.i.d. unit Fr\'echet
  random variables which are also chosen to be independent of $X$ and
  let $\xi$ be the coordinatewise maximum of $c_u u\zeta_u$ over all
  $u\in M$. It remains to show that $\xi$ has the required
  distribution. Consider an arbitrary point $v\in M$. By the condition
  on $M$, $tu\leq v$ is possible for some $t>0$ if and only if $u\leq
  v$ and $t\leq 1$. Thus,
  \begin{align*}
    \Prob{\xi\leq v}&=\prod_{u\in M} \Prob{c_u\zeta_u u\leq v}
    =\prod_{u\in M,\; u\leq v}\Prob{c_u\zeta_u\leq1}
    =\exp\left\{-\sum_{u\in M,\; u\leq v} c_u\right\}\\
    &=\exp\left\{-h(\vee M)\Prob{X\cap\{u:\; u\leq v\}\neq\emptyset}\right\}\\
    &=\exp\left\{-h(\vee M)T(\{u:\; u\leq v\})\right\}
    =\exp\left\{-h(v)\right\}\,.
  \end{align*}
\end{proof}

A simple example of set $M$ from Theorem~\ref{thr:consist} is the
smallest set which contains all basis vectors in $\R^d$ and is closed
with respect to coordinatewise maxima. Then $M$ consists of the
vertices of the unit cube $\II^d$ without the origin and the values of
$h$ on $M$ become the extremal coefficients.  The \emph{extremal
  coefficients} $\theta_A$ of a simple max-stable random vector
$\xi=(\xi_1,\dots,\xi_d)$ are defined from the equations
\begin{equation}
  \label{eq:theta-A}
  \Prob{\max_{j\in A} \xi_j\leq z}=(\Prob{\xi_1\leq z})^{\theta_A}\,,
  \quad z>0\,,\; A\subset\{1,\dots,d\}\,,
\end{equation}
see \cite{sch:taw02,sch:taw03}. Since the marginals are unit
Fr\'echet, it suffices to use (\ref{eq:theta-A}) for $z=1$ only. If
$e_A=\sum_{i\in A} e_i$, then (\ref{eq:fx=e-x-quad}) implies 
\begin{displaymath}
  \theta_A=h(K,e_A)=\|e_A\|_K\,.
\end{displaymath}
Every nonempty set $A\subset\{1,\dots,d\}$ can be associated with the
unique vertex of the unit cube $\II^d\setminus\{0\}$. The consistency
condition for the extremal coefficients follows directly from
Theorem~\ref{thr:consist} and can be formulated as follows.

\begin{cor}
  \label{cor:thetaa}
  A family of non-negative numbers $\theta_A$,
  $A\subset\{1,\dots,d\}$, is a set of extremal coefficients for a
  simple max-stable distribution if and only if $\theta_{\emptyset}=0$
  and $\theta_A$ is a union-completely alternating function of $A$.
\end{cor}

This consistency result for the extremal coefficients has been
formulated in \cite{sch:taw02} as a set of inequalities that, in fact,
mean the complete alternation property of $\theta_A$.

\section{Operations with dependency sets}
\label{sec:oper-with-conv}

\paragraph{Rescaling.} For a dependency set $K$ and
$\lambda_1,\dots,\lambda_d>0$ define
\begin{equation}
  \label{eq:lambda-k=lambd-dots}
  \lambda K=\{(\lambda_1x_1,\dots,\lambda_dx_d):\;
  x=(x_1,\dots,x_d)\in K\}\,.
\end{equation}
Then $e^{-h(\lambda K,x^*)}$ is the cumulative distribution function
of $\lambda^*\xi=(\xi_1/\lambda_1,\dots,\xi_d/\lambda_d)$. 

\paragraph{Projection.}
If $\xi'$ denotes the vector composed from the first $k$-coordinates
of $d$-dimensional vector $\xi$ with the dependency set $K$, then
\begin{align*}
  \Prob{\xi'\leq (x_1,\dots,x_k)}
  &=\exp\{-\|(x_1,\dots,x_k,\infty,\dots,\infty)^*\|_K\}\\
  &=\exp\{-\|(x_1,\dots,x_k)^*\|_{K'}\}\,,
\end{align*}
where $K'$ is the projection of $K$ onto the subspace spanned by the
first $k$ coordinates in $\R^d$. Thus, taking a sub-vector of $\xi$
corresponds to projecting of $K$ onto the corresponding coordinate
subset. Recall Proposition~\ref{prop:mzp} which says that if all
two-dimensional projections of $K$ are squares, then $K$ is
necessarily the cube.

\begin{prop}
  \label{prop:proj-sect}
  If\; $\LL$ is the subspace spanned by some coordinate axes in $\R^d$,
  the projection of $K$ onto $\LL$ coincides with $K\cap \LL$.
\end{prop}
\begin{proof}
  By definition, $K=c\E\simplex{\eta}$. Then it suffices to note that
  the projection of $\simplex{\eta}$ on $\LL$ equals
  $\simplex{\eta}\cap \LL$. Indeed every selection of
  $\simplex{\eta}\cap\LL$ can be associated with the projection of a
  selection of $\simplex{\eta}$.
\end{proof}

An interesting open question concerns a reconstruction of $K$ from its
lower-dimensional projections. In various forms this question was
discussed in \cite[Sec.~3.5.6]{kot:nad00} and \cite[Sec.~4.7]{joe97}.

\paragraph{Cartesian product.} 
If $K'$ and $K''$ are two dependency sets of simple max-stable random
vectors $\xi'$ and $\xi''$ with dimensions $d'$ and $d''$
respectively, then the Cartesian product $K'\times K''$ is the
dependency set corresponding to the max-stable random vector $\xi$
obtained by concatenating of independent copies of $\xi'$ and $\xi''$.
Indeed, if $x=(x',x'')$, then
\begin{align*}
  \Prob{\xi\leq x} &= \exp\{-h(K'\times K'',x)\}
  =\exp\{-h(K',x')-h(K'',x'')\}\\
  &=\Prob{\xi'\leq x'}\Prob{\xi''\leq x''}.
\end{align*}

\paragraph{Minkowski sum.} 
If $K'$ and $K''$ are dependency sets of two independent max-stable
random vectors $\xi'$ and $\xi''$ of dimension $d$, then the weighted
Minkowski sum $K=\lambda K'+(1-\lambda)K''$ with $\lambda\in[0,1]$ is
the dependency set of the max-stable random vector
\begin{equation}
  \label{eq:xi-rep-m}
  \xi=(\lambda\xi')\vee((1-\lambda)\xi'')\,.
\end{equation}
The cumulative distribution functions of $\xi',\xi''$ and $\xi$ are
related as
\begin{displaymath}
  F_\xi(x)=F_{\xi'}(x)^\lambda F_{\xi''}(x)^{(1-\lambda)}\,.
\end{displaymath}

It is possible to generalise the Minkowski summation scheme for
multivariate weights.
Consider $K=\lambda K'+(1-\lambda)K''$ for some $\lambda\in[0,1]^d$,
where the products of vectors and sets are defined in
(\ref{eq:lambda-k=lambd-dots}).  Then
$\|x\|_K=\|\lambda^*x\|_{K'}+\|(1-\lambda)^*x\|_{K''}$, so that
(\ref{eq:xi-rep-m}) also holds with the products defined
coordinatewisely.

\begin{example}
  \label{ex:2d-cut}
  If $\xi_1$ and $\xi_2$ are independent with unit Fr\'echet
  distributions and $\alpha_1,\alpha_2\in[0,1]$, then setting
  $\lambda=(\alpha_1,1-\alpha_2)$ we obtain the max-stable random vector
  \begin{displaymath}
    \xi=(\alpha_1\xi_1\vee(1-\alpha_1)\xi_2
    ,(1-\alpha_2)\xi_1\vee\alpha_2\xi_2)
  \end{displaymath}
  with the dependency set
  $K=\conv\{(0,0),(0,1),(1,0),(\alpha_1,1),(1,\alpha_2)\}$. If
  $\alpha_1=\alpha_2=\alpha$, then $\xi$ has the Marshall-Olkin
  distribution, cf \cite[Ex.~4.1.1]{fal:hus:rei04}.
\end{example}

\begin{example}[Matrix weights]
  \label{ex:mw}
  Let $a_{ij}$, $i=1,\dots,m$, $j=1,\dots,d$, be a matrix of positive
  numbers such that $\sum_{i=1}^m a_{ij}=1$ for all $j$. Furthermore,
  let $\zeta_1,\dots,\zeta_m$ be i.i.d. random variables with
  $\Phi_1$-distribution. Define $\xi=(\xi_1,\dots,\xi_d)$ by
  \begin{displaymath}
    \xi_j=\max_{1\leq i\leq m} \zeta_i a_{ij}\,, \quad j=1,\dots,d\,,
  \end{displaymath}
  cf \cite[Lemma~4.1.2]{fal:hus:rei04}.  Then $\xi$ is simple
  max-stable with the corresponding norm 
  \begin{displaymath}
    \|x\|_K=\sum_{i=1}^m \max_{1\leq j\leq d} a_{ij}x_j\,,
  \end{displaymath}
  i.e. its dependency set is 
  \begin{math}
    K=\simplex{(a_{11},\dots,a_{1d})}
    +\cdots+\simplex{(a_{m1},\dots,a_{md})}\,. 
  \end{math}
\end{example}

\paragraph{Power sums.}
A \emph{power-mean} of two convex compact sets $K'$ and $K''$
containing the origin in their interior is defined to be a convex set
$K$ such that
\begin{equation}
  \label{eq:hk-up=lambda-hk}
  h(K,x)^p=\lambda h(K',x)^p+(1-\lambda)h(K,x)^p\,, 
\end{equation}
where $\lambda\in[0,1]$ and $p\geq1$, see \cite{fir67}. 
The power-mean definition is applicable also if $K'$ and $K''$ satisfy
(\ref{eq:ksub}), despite the fact that the origin is not their
interior point. In the plane, the power sum is a dependency set if
$K'$ and $K''$ satisfy (\ref{eq:ksub}). Therefore, the power sum of
dependency sets leads to a new operation with distributions of
\emph{bivariate} max-stable random vectors. For instance, if $K'$ is
the unit cross-polytope and $K''$ is the unit square, then, for $p=2$,
\begin{displaymath}
  \|x\|_K=((x_1+x_2)^2+(\max(x_1,x_2))^2)^{1/2}\,.
\end{displaymath}

\paragraph{Minkowski difference.} Let $K'$ and $K''$ be two
dependency sets. For any $\lambda>0$ define
\begin{displaymath}
  L=K'-\lambda K''
  =\{x:\; x+\lambda K''\subset K'\}\,.
\end{displaymath}
If the spectral measures $\sigma'$ and $\sigma''$ of $K'$ and $K''$
are such that $\sigma=\sigma'-\lambda\sigma''$ is a non-negative
measure, then $L$ is a max-zonoid regardless of the dimension of the
space. The negative logistic distribution from Example~\ref{ex:neg-lo}
illustrates this construction.

\paragraph{Convex hull and intersection.}
In the space of a dimension $d\geq3$ the convex hull or intersection
of dependency sets do not necessarily remain dependency sets.
However, on the plane this is always the case.

Let $K'$ and $K''$ be the dependency sets of bivariate simple
max-stable random vectors $\xi'$ and $\xi''$. Since $h(\conv(K'\cup
K''),x)=h(K',x)\vee h(K'',x)$, the dependency set $K=\conv(K'\cup
K'')$ corresponds to a max-stable random vector $\xi$ such that
\begin{displaymath}
  \Prob{\xi\leq x}=\min(\Prob{\xi'\leq x},\Prob{\xi''\leq x})\,,
  \quad x\in[0,\infty)^2\,. 
\end{displaymath}

The intersection of two planar dependency sets also remains the
dependency set and so yields another new operation with distributions
of simple max-stable bivariate random vectors.

\paragraph{Duality.}
If the polar to the dependency set $K$ of $\xi$ is a max-zonoid, then
the corresponding simple max-stable random vector $\xi^o$ is said to
be the dual to $\xi$.  In the plane, the polar to a max-zonoid is
max-zonoid; it is not known when it holds in higher dimensions. This
duality operation is a new operation with distributions of bivariate
max-stable random vectors, see also Example~\ref{ex:basic}.

\section{Infinite dimensional case}
\label{sec:infin-dimens-case}

It is possible to define the dependency set for max-stable stochastic
processes studied in \cite{haa84,fal:hus:rei04,gin:hah:vat90}. The
spectral representation \cite[Prop.~3.2]{gin:hah:vat90} of a sample
continuous max-stable process $\xi(t)$, $t\in S$, on a compact metric
space $S$ with unit Fr\'echet marginals yields that
\begin{displaymath}
  -\log \Prob{\xi< f}=\int_{\SS} \|g/f\|_\infty d\sigma(g)\,,
\end{displaymath}
where $\SS$ is the family of non-negative continuous functions $g$ on
$S$ that their maximum value $\|g\|_\infty$ equals $1$, and $\sigma$
is a finite Borel measure on $\SS$ such that $\int_\SS gd\sigma(g)$ is
the function identically equal to $1$.

The corresponding dependency set is the set in the space of finite
measures with the total variation distance, which is the dual space to
the family of non-negative continuous functions.  For a continuous
function $g$, define $\simplex{g}$ to be the closed convex hull of the
family of atomic measures $g(x)\delta_x$ for $x\in S$.  Then the
dependency set is the expectation of $c\simplex{\eta}$, where
$c=\sigma(\SS)$ and $\eta$ is distributed according to the normalised
$\sigma$.

\section*{Acknowledgements}
\label{sec:acknowledgements}

The problematic of this paper was motivated by a lecture of
Prof.~M.~Falk given in Bern, see
\textsf{http://statistik.mathematik.uni-wuerzburg.de/$\sim$falk/archiv/bern\_04.pdf}.
The author is grateful to J\"urg H\"usler for discussions on this
topic and extreme values in general. The referee's comments led
to clarification of arguments as well as corrections and
simplification of several proofs.

This work was supported by the Swiss National Science Foundation Grant
No.  200021-103579.

\newcommand{\noopsort}[1]{} \newcommand{\printfirst}[2]{#1}
  \newcommand{\singleletter}[1]{#1} \newcommand{\switchargs}[2]{#2#1}

\end{document}